\documentclass[a4paper,reqno,12pt]{amsart}
\usepackage{enumerate}
\usepackage{amssymb}
\usepackage{ifthen}
\usepackage{graphicx}
\usepackage{color,xcolor}
\usepackage{mathrsfs}
\numberwithin{equation}{section}
\newtheorem{thm}{Theorem}
\newtheorem{theorem}{Theorem}
\newtheorem{lem}{Lemma}
\newtheorem{cor}{Corollary}
\newtheorem{defn}{Definition}
\newtheorem{example}{Example}

\newtheorem{claim}{Claim}

\newtheorem{conj}{Conjecture}

\newtheorem{prob}{Problem}

\newenvironment{rem}{%
\bigskip
\noindent \textsl{{\sl Remark. }}}{\bigskip}
\newenvironment{rems}{%
\bigskip
\noindent \textsl{{\sl Remarks. }}}{\bigskip}

\newenvironment{pf}[1][]{%
 \vskip 1mm
 \noindent
 \ifthenelse{\equal{#1}{}}%
  {{\slshape Proof. }}%
  {{\slshape #1.} }%
 }%
{\qed\bigskip}

\newcounter{alphabet}
\newcounter{tmp}
\newenvironment{Thm}[1][]{\refstepcounter{alphabet}%
\bigskip%
\noindent%
{\bf Theorem \Alph{alphabet}}%
\ifthenelse{\equal{#1}{}}{}{ (#1)}%
{\bf .} \itshape}{\vskip 8pt}

\makeatletter
\newcommand{\Ref}[1]{\@ifundefined{r@#1}{}{\setcounter{tmp}{\ref{#1}}\Alph{tmp}}}
\makeatother

\newcommand{\IN}{{\mathbb N}}
\newcommand{\IC}{{\mathbb C}}
\newcommand{\ID}{{\mathbb D}}

\def\be{\begin{equation}}
\def\ee{\end{equation}}
\def\bes{\begin{equation*}}
\def\ees{\end{equation*}}

\newcommand{\bee}{\begin{enumerate}}
\newcommand{\eee}{\end{enumerate}}

\newcommand{\blem}{\begin{lem}}
\newcommand{\elem}{\end{lem}}
\newcommand{\bthm}{\begin{thm}}
\newcommand{\ethm}{\end{thm}}
\newcommand{\bcor}{\begin{cor}}
\newcommand{\ecor}{\end{cor}}
\newcommand{\beg}{\begin{example}}
\newcommand{\eeg}{\end{example}}
\newcommand{\begs}{\begin{examples}}
\newcommand{\eegs}{\end{examples}}
\newcommand{\bdefe}{\begin{defn}}
\newcommand{\edefe}{\end{defn}}
\newcommand{\bprob}{\begin{prob}}
\newcommand{\eprob}{\end{prob}}
\newcommand{\bques}{\begin{ques}}
\newcommand{\eques}{\end{ques}}
\newcommand{\bei}{\begin{itemize}}
\newcommand{\eei}{\end{itemize}}

\newcommand{\bde}{\begin{deter}}
\newcommand{\ede}{\end{deter}}
\newcommand{\bca}{\begin{case}}
\newcommand{\eca}{\end{case}}
\newcommand{\bcl}{\begin{claim}}
\newcommand{\ecl}{\end{claim}}
\newcommand{\bcon}{\begin{conj}}
\newcommand{\econ}{\end{conj}}
\newcommand{\bcons}{\begin{conjs}}
\newcommand{\econs}{\end{conjs}}
\newcommand{\bprop}{\begin{propo}}
\newcommand{\eprop}{\end{propo}}
\newcommand{\br}{\begin{rem}}
\newcommand{\er}{\end{rem}}
\newcommand{\brs}{\begin{rems}}
\newcommand{\ers}{\end{rems}}
\newcommand{\bo}{\begin{obser}}
\newcommand{\eo}{\end{obser}}
\newcommand{\bos}{\begin{obsers}}
\newcommand{\eos}{\end{obsers}}
\newcommand{\bpf}{\begin{pf}}
\newcommand{\epf}{\end{pf}}
\newcommand{\ba}{\begin{array}}
\newcommand{\ea}{\end{array}}
\newcommand{\beq}{\begin{eqnarray}}
\newcommand{\beqq}{\begin{eqnarray*}}
\newcommand{\eeq}{\end{eqnarray}}
\newcommand{\eeqq}{\end{eqnarray*}}

\newcommand{\ds}{\displaystyle}

\newcounter{minutes}
\divide\time by 60
\newcounter{hours}
\multiply\time by 60 \addtocounter{minutes}{-\time}

\begin{document}
\bibliographystyle{amsplain}
\title[Refined Bohr inequality for bounded analytic functions]
{Refined Bohr inequality for bounded analytic functions}

\def\thefootnote{}
\footnotetext{ \texttt{\tiny File:~\jobname .tex,
          printed: \number\day-\number\month-\number\year,
          \thehours.\ifnum\theminutes<10{0}\fi\theminutes}
} \makeatletter\def\thefootnote{\@arabic\c@footnote}\makeatother

\author[G. Liu]{Gang Liu}
\address{G. Liu, College of Mathematics and Statistics
 (Hunan Provincial Key Laboratory of Intelligent Information Processing and Application),
Hengyang Normal University, Hengyang,  Hunan 421002, China}
\email{liugangmath@sina.cn}

\author[Z.H. Liu]{Zhihong Liu}
\address{Z.H. Liu, College of Science, Guilin University of Technology, Guilin, Guangxi 541004, China}
\email{liuzhihongmath@163.com}

\author[S. Ponnusamy]{Saminathan Ponnusamy
}

\address{S. Ponnusamy, Department of Mathematics,
Indian Institute of Technology Madras, Chennai-600 036, India}
\email{samy@iitm.ac.in}



\begin{abstract}
In this article, by combining appropriate refined Bohr's inequalities with some techniques concerning bounded analytic functions defined in the unit disk,
we generalize and improve several Bohr type inequalities for such functions.
\end{abstract}

\maketitle \pagestyle{myheadings}
\markboth{G. Liu, Z. H. Liu and S. Ponnusamy}{Refined Bohr inequality for bounded analytic functions}

\section{Introduction}  \label{sec1}
Throughout the paper, ${\mathcal B}$  denote the set of all analytic self-maps of the unit disk  $\mathbb{D}=\{z\in \IC:\, |z| < 1 \}$.
Let us start recalling the remarkable discovery of H. Bohr in 1914 (see \cite{boh}).

\begin{Thm}\label{thmA}
If $f\in {\mathcal B}$ and $f(z)=\sum_{n=0}^{\infty} a_n z^n$, then the following sharp inequality holds:
$ \sum_{n=0}^{\infty} |a_n|\, r^n \leq 1$ for  $r\leq 1/3$.
\end{Thm}

 We remark that if $|f(z)|\leq 1$ in $\ID$ and $|f(z_0)|=1$ for some point $z_0\in\ID$, then  $f(z)$ reduces a unimodular constant
and thus,  we restrict our attention to $f\in{\mathcal B}$.
The above inequality and the sharp constant $1/3$ are known as the classical Bohr inequality and the Bohr radius for the family ${\mathcal B}$, respectively.
Bohr \cite{boh}  showed it only for $r\leq1/6$ and this article includes Wiener's proof showing that $r=1/3$ is sharp.
Subsequently, some different proofs were given (see \cite{PPS,PS2004,PS2006,sid,tom} and also the recent survey chapters \cite{AAP} and \cite[Chapter 8]{GarMasRoss-2018}).
It is worth pointing out that if $|a_0|$ in Bohr inequality is replaced by $|a_0|^2$,
then the constant $1/3$ could be replaced by $1/2$.
Moreover, if $a_0=0$ in Theorem \Ref{thmA}, then the sharp Bohr radius can be improved to be $1/\sqrt{2}$
(see for example \cite{KP2017}, \cite[Corollary~2.9]{PPS} and the recent paper \cite{PW} for a general result).
For these results, one of the proofs relied on the sharp coefficient inequalities, i.e., $|a_n|\leq1-|a_0|^2$ $(n\geq1, f\in{\mathcal B})$.
By means of these inequalities, it is pointed out in \cite{KP2017} that the sharp result in Theorem \Ref{thmA} cannot be obtained in the
extremal case $|a_0|<1$.
However, the sharp version of Theorem \Ref{thmA} has been achieved for any individual function from ${\mathcal B}$ (see \cite{AKP})
and  for some  subclasses of univalent functions, we refer to \cite{abu,aiz}. In \cite{PVW} (see also \cite{PVW2}), the authors established the following
refined Bohr inequality by applying a refined version of the coefficient inequalities
carefully (see Lemma \ref{lem3}).

\begin{Thm}{ \rm(\cite[Theorem~1]{PVW})} \label{thmB}
Suppose that $f\in {\mathcal B}$ and $f(z)=\sum_{n=0}^{\infty} a_n z^n$.
Then
\begin{equation*}
 \sum_{n=0}^\infty |a_n|r^n+  \left(\frac{1}{1+|a_0|}+\frac{r}{1-r}\right)\sum_{n=1}^{\infty}|a_n|^2r^{2n} \leq 1 ~\mbox{for}~ r \leq \frac{1}{2+|a_0|}
\end{equation*}
and the numbers $\frac{1}{2+|a_0|}$ and $\frac{1}{1+|a_0|}$ cannot be improved.  Moreover,
\begin{equation*}
|a_0|^2+ \sum_{n=1}^\infty |a_n|r^n+ \left(\frac{1}{1+|a_0|}+\frac{r}{1-r}\right)\sum_{n=1}^{\infty}|a_n|^2r^{2n}\leq 1 ~\mbox{for}~ r \leq \frac{1}{2}
\end{equation*}
and the numbers $\frac{1}{2}$ and $\frac{1}{1+|a_0|}$ cannot be improved.
\end{Thm}

Besides these results, there are plenty of works about Bohr inequality for the family ${\mathcal B}$.
Based on the notion of Rogosinski's inequality and Rogosinski's radius investigated in \cite{LG,rog,SS}, Kayumov and Ponnusamy \cite{KP-pre} introduced and obtained
the following Bohr-Rogosinski inequality and Bohr-Rogosinski radius.

\begin{Thm}{ \rm(\cite[Theorem~1]{KP-pre})}\label{thmC}
Suppose that $f\in {\mathcal B}$ and $f(z)=\sum_{n=0}^{\infty} a_n z^n$. Then
$$
|f(z)|+\sum_{n=N}^{\infty}\left|a_{n}\right| r^{n} \leq 1 \quad \text {for} ~|z|=r\leq R_{N},
$$
where $R_{N}$ is the positive root of the equation  $2(1+r) r^{N}-(1-r)^{2}=0$.
The radius $R_{N}$ is best possible. Moreover,
$$
|f(z)|^{2}+\sum_{n=N}^{\infty}\left|a_{n}\right| r^{n} \leq 1 \quad \text {for} ~|z|=r\leq R_{N}^{\prime},
$$
where $R_{N}^{\prime}$ is the positive root of the equation $(1+r) r^{N}-(1-r)^{2}=0 .$ The radius $R_{N}^{\prime}$ is
best possible.
\end{Thm}

In order to determine the Bohr radius for the class of odd functions in the family ${\mathcal B}$, which was posed in \cite{ABS},
the authors in \cite{KP2017,KP2018-1} established a more general result consisting of functions of the form $f_{p,m}(z)=\sum_{n=0}^{\infty} a_{pn+m} z^{pn+m}$ in ${\mathcal B}$.
Note that $f_{p,1}$ is a $p-$symmetric function and $f_{2,1}$ is an odd function.

\begin{Thm}{ \rm(\cite[Theorem~1]{KP2018-1})}\label{thmD}
Let p $\in \mathbb{N}$ and 0 $\leq m \leq p$.
Suppose that $f\in {\mathcal B}$ and $f(z)=\sum_{n=0}^{\infty} a_{pn+m} z^{pn+m}$. Then
$$ \sum_{n=0}^{\infty}\left|a_{pn+m}\right| r^{pn+m} \leq 1 \quad \text {for} ~r \leq r_{p, m},
$$
where $r_{p, m}$ is the maximal positive root of the equation
$-6 r^{p-m}+r^{2(p-m)}+8 r^{2 p}+1=0.
$
The extremal function has the form $z^{m}\left(z^{p}-a\right) /\left(1-az^{p}\right),$ where
$$a=\left(1-\frac{\sqrt{1-r_{p, m}^{2 p}}}{\sqrt{2}}\right) \frac{1}{r_{p, m}^{p}}.
$$
\end{Thm}

 Besides these results,  few different formulations of improved Bohr inequalities were
obtained in \cite{KP2018-2} and also in the recent articles \cite{PVW,PW}.   Below we recall one of them which is associated
with area.

\begin{Thm}{ \rm(\cite[Theorem~1]{KP2018-2})} \label{thmE}
Suppose that $f\in {\mathcal B}$, $f(z)=\sum_{n=0}^{\infty} a_n z^n$ and $S_r$ denotes the area of the
Riemann surface of the function $f^{-1}$ defined on the image of the subdisk $|z|<r$ under the mapping $f$. Then
$$\sum_{n=0}^{\infty}\left|a_{n}\right| r^{n}+\frac{16}{9}\left(\frac{S_{r}}{\pi}\right) \leq 1~\text{for}~ r\leq\frac{1}{3}
$$
and the numbers $1/3$ and $16/9$ cannot be improved. Moreover,
$$
\left|a_{0}\right|^{2}+\sum_{n=1}^{\infty}\left|a_{n}\right| r^{n}+\frac{9}{8}\left(\frac{S_{r}}{\pi}\right) \leq 1 \quad \text {for} ~r\leq\frac{1}{2}
$$
and the constants $1/2$ and $9/8$ cannot be improved.
\end{Thm}

 Based on the initiation of Kayumov and Ponnusamy\cite{KP2018-2}, several forms of Bohr-type inequalities for the family ${\mathcal B}$
were considered in \cite{LSX} when the  Taylor coefficients of classical Bohr inequality are partly or completely replaced by higher order derivatives of $f$.
Here we only recall one of them.

\begin{Thm}{ \rm(\cite[Theorem~2.1]{LSX})} \label{thmF}
Suppose that $f\in {\mathcal B}$ and $f(z)=\sum_{n=0}^{\infty} a_n z^n$. Then the following sharp inequality holds:
$$
|f(z)|+\left|f^{\prime}(z)\right|r+\sum_{n=2}^{\infty}\left|a_{n}\right|r^{n} \leq 1 \quad \text {for}~|z|=r \leq \frac{\sqrt{17}-3}{4}.
$$
\end{Thm}

Besides these several authors have investigated some other extensions along with applications and connections with local Banach space theory and other topics.
See for instance \cite{Aizen-00-1,AAD,BDK5,BombBour-196,BombBour-2004,DF,DGM,DeGarM04,DMS,DP,GarMasRoss-2018} and the references therein.
In view of the recent developments on Bohr-type inequalities, it is natural to ask whether
we can further generalize or improve these results in the recent setting of \cite{PVW}. In this article, we give
an affirmative answer to this question.

The paper is organized as follows. In Section \ref{sec2}, we include some
preliminary results and a new lemma. In Section \ref{sec3}, we establish an improved version of the Bohr-Rogosinski inequality.
In Section \ref{sec4}, we obtain refined Bohr inequalities for symmetric functions.
Finally, in Section \ref{sec5}, we consider some refined Bohr type inequalities associated with
area, modulus of $f-a_0(f)$ and higher order derivatives of $f$ in part.

\section{Preliminary results}  \label{sec2}

In order to obtain our results, we need few lemmas.

\begin{lem} { \rm(\cite[Proof~ of~ Theorem ~1]{KP2017} and \cite{KP2018-1})} \label{lem1}
 Suppose that $f\in {\mathcal B}$ and $f(z)=\sum_{n=0}^{\infty} a_n z^n$. Then we have
$$ \sum_{n=1}^{\infty}\left|a_{n}\right| r^{n} \leq
\begin{cases}
\displaystyle r\frac{1-|a_0|^2}{1-r|a_0|}, &\, \text{for}\quad |a_0|\geq r,\\
 \displaystyle r\frac{\sqrt{1-|a_0|^2}}{\sqrt{1-r^2}}, &\, \text{for}\quad |a_0|<r.
\end{cases}
$$
\end{lem}

\begin{lem} { \rm(\cite[Lemma~1]{KP2018-2})} \label{lem2}
Suppose that $f\in {\mathcal B}$ and $f(z)=\sum_{n=0}^{\infty} a_n z^n$. Then the
following sharp inequality holds:
$$\sum_{n=1}^{\infty} n\left|a_{n}\right|^{2} r^{2n} \leq r^{2} \frac{\left(1-\left|a_{0}\right|^{2}\right)^{2}}{\left(1-\left|a_{0}\right|^{2} r^{2}\right)^{2}}
~\mbox{ for  $\ds 0<r\leq\frac{1}{\sqrt{2}}.$}
$$
\end{lem}

\begin{lem}  { \rm(\cite{Car} and \cite[Lemma~B]{PVW})} \label{lem3}
Suppose that $f\in {\mathcal B}$ and $f(z)=\sum_{n=0}^{\infty} a_n z^n$. Then the following inequalities hold.
\begin{enumerate}
	\item[{\rm (a)}] $|a_{2n+1}|\leq 1-|a_0|^2-\cdots - |a_n|^2,~\ n=0,1,\ldots$ 
	\item[{\rm (b)}] $|a_{2n}|\leq 1-|a_0|^2-\cdots -|a_{n-1}|^2 - \frac{|a_n|^2}{1+|a_0|} ,~\ n=1,2,\ldots$. 
\end{enumerate}
Further, to have equality in {\rm (a)} it is necessary that $f$ is a rational function of the form
$$ f(z)=\frac{a_{0}+a_{1}z+ \cdots + a_{n}z^{n}+\epsilon z^{2n+1}}{1+(\overline{a_n}z^n+ \cdots +\overline{a_0}z^{2n+1})\epsilon},~|\epsilon|=1,
$$
and  to have equality in {\rm (b)} it is necessary that $f$ is a rational function of the form
$$ f(z)=\frac{a_{0}+a_{1}z+ \cdots + \frac{a_{n}}{1+|a_0|} z^{n}+ \epsilon z^{2n} }{1+\left(\frac{\overline{a_n}}{1+|a_0|}z^n+ \cdots +\overline{a_0}z^{2n}\right) \epsilon},
~|\epsilon|=1,
$$
where  $a_0 \overline{a_n}^2 \epsilon$ is non-positive real.
\end{lem}

In what follows, $\lfloor x\rfloor $ denotes the largest integer no more than $x$, where $x$ is a real number.

\begin{lem}\label{lem4}
Suppose that $f\in {\mathcal B}$ and $f(z)=\sum_{n=0}^{\infty} a_n z^n$.
Then for any $N\in\mathbb{N}$, the following inequality holds:
$$	\sum_{n=N}^{\infty}|a_n|r^n+sgn(t)\sum_{n=1}^{t}|a_n|^2\frac{r^N}{1-r}+ \left(\frac{1}{1+|a_0|}
+\frac{r}{1-r}\right)\sum_{n=t+1}^{\infty} |a_n|^2 r^{2n}\leq(1-|a_0|^2)\frac{r^N}{1-r},
$$
for $r\in[0,1)$, where $t=\lfloor(N-1)/2\rfloor$.
\end{lem}

\bpf
The proof is divided into two cases.
For the case of   even values of $N\in \IN$, we set $N=2m$ ($m\geq 1$).  It follows from Lemma \ref{lem3} that
\begin{align*}
\sum_{n=2m}^{\infty}|a_n|r^n=  \,&  \sum_{n=m}^{\infty}|a_{2n}|r^{2n}+\sum_{n=m}^{\infty}|a_{2n+1}|r^{2n+1}\\
\leq \, & \sum_{n=m}^{\infty}\left(1-\sum_{k=0}^{n-1}|a_k|^2 - \frac{|a_n|^2}{1+|a_0|}\right)r^{2n}
+\sum_{n=m}^{\infty}\left(1-\sum_{k=0}^n|a_k|^2\right)r^{2n+1}\\
=\, &\sum_{n=2m}^{\infty}r^n-\left(\sum_{n=0}^{m-1}|a_n|^2\right)\left(\sum_{n=2m}^{\infty}r^n\right)-\left(\frac{1}{1+|a_0|}
+\sum_{n=1}^{\infty}r^n\right)\sum_{n=m}^{\infty}|a_n|^2r^{2n}\\
=\, &\left(1-\sum_{n=0}^{m-1}|a_n|^2\right)\frac{r^{2m}}{1-r}-\left(\frac{1}{1+|a_0|}
+\frac{r}{1-r}\right)\sum_{n=m}^{\infty} |a_n|^2 r^{2n}
\end{align*}
so that
$$ \sum_{n=2m}^{\infty}|a_n|r^n  +\left (\sum_{n=1}^{m-1}|a_n|^2\right)\frac{r^{2m}}{1-r}+\left(\frac{1}{1+|a_0|}
+\frac{r}{1-r}\right)\sum_{n=m}^{\infty} |a_n|^2 r^{2n} \leq(1-|a_0|^2)\frac{r^N}{1-r},
$$
and the desired result follows easily.

For the case of  odd values of $N\in \IN$  associated with the  inequality in the statement,  we set $N=2m+1$ ($m\geq 0$).
It follows from Lemma \ref{lem3} that
\begin{align*}
\sum_{n=2m+1}^{\infty}|a_n|r^n=\,
& \sum_{n=2(m+1)}^{\infty}|a_{n}|r^{n}+|a_{2m+1}|r^{2m+1}\\
\leq \, & \left(1-\sum_{n=0}^{m}|a_n|^2\right)\frac{r^{2(m+1)}}{1-r}-\left(\frac{1}{1+|a_0|}
+\frac{r}{1-r}\right)\sum_{n=m+1}^{\infty} |a_n|^2 r^{2n}\\
& \hspace{1cm} +\left(1-\sum_{n=0}^m|a_n|^2\right)r^{2m+1}\\
=\, &\left(1-\sum_{n=0}^{m}|a_n|^2\right)\frac{r^{2m+1}}{1-r}-\left(\frac{1}{1+|a_0|}
+\frac{r}{1-r}\right)\sum_{n=m+1}^{\infty} |a_n|^2 r^{2n}.
\end{align*}
Simple translation gives the desired result. This completes the proof.
\epf

\section{Refined Bohr-Rogosinski inequalities}  \label{sec3}

Using Lemma \ref{lem4} and the similar proof of Theorem \Ref{thmC}, we can easily get the following result.
So we omit the details.

\begin{theorem} \label{thm1}
Suppose that $f\in {\mathcal B}$ and $f(z)=\sum_{n=0}^{\infty} a_n z^n$.  For $n\geq \IN$, let $t=\lfloor(N-1)/2\rfloor$.
Then $$|f(z)|+\sum_{n=N}^{\infty}|a_n|r^n+sgn(t)\sum_{n=1}^{t}|a_n|^2\frac{r^N}{1-r}+ \left(\frac{1}{1+|a_0|}
+\frac{r}{1-r}\right)\sum_{n=t+1}^{\infty} |a_n|^2 r^{2n} \leq 1
$$
for $|z|=r\leq R_{N}$, where $R_{N}$ is    as in Theorem {\rm \Ref{thmC}}.
The radius $R_{N}$ is best possible. Moreover,
$$ |f(z)|^{2}+\sum_{n=N}^{\infty}|a_n|r^n+sgn(t)\sum_{n=1}^{t}|a_n|^2\frac{r^N}{1-r}+ \left(\frac{1}{1+|a_0|}
+\frac{r}{1-r}\right)\sum_{n=t+1}^{\infty} |a_n|^2 r^{2n} \leq 1
$$
for $|z|=r\leq R_{N}^{\prime}$, where $R_{N}^{\prime}$ is  as in Theorem {\rm \Ref{thmC}}.
The radius $R_{N}^{\prime}$ is best possible.
\end{theorem}

For the case of $N=1$, it is easy to see that $R_1=\sqrt{5}-2$ and $R_1^{'}=1/3$.
However, the two constants can be improved  for any individual function in ${\mathcal B}$ (in the context of Theorem \Ref{thmC}
and Theorem \ref{thm1} with $N=1$).

\begin{theorem} \label{thm2}
Suppose that $f\in {\mathcal B}$ and $f(z)=\sum_{n=0}^{\infty} a_n z^n$. Then 
$$A(z):=|f(z)|+\sum_{n=1}^{\infty}|a_n|r^n+ \left(\frac{1}{1+|a_0|}
+\frac{r}{1-r}\right)\sum_{n=1}^{\infty} |a_n|^2 r^{2n} \leq 1
$$
for $|z|=r\leq r_{a_0}= 2/\left(3+|a_0|+\sqrt{5}(1+|a_0|)\right).$
The radius $r_{a_0}$ is best possible and $r_{a_0}\geq \sqrt{5}-2$. Moreover,
$$B(z):= |f(z)|^2+\sum_{n=1}^{\infty}|a_n|r^n+ \left(\frac{1}{1+|a_0|}
+\frac{r}{1-r}\right)\sum_{n=1}^{\infty} |a_n|^2 r^{2n} \leq 1
$$
for $|z|=r\leq r_{a_0}'$,
where $r_{a_0}'$ is the unique positive root of the equation
$$(1-|a_0|^3)r^3-(1+2|a_0|)r^2-2r+1=0.
$$
The radius $r_{a_0}'$ is best possible.
Further, we have $1/3<r_{a_0}'<1/(2+|a_0|)$. 
\end{theorem}

\bpf We consider the first part.
According to Lemma \ref{lem4} with $N=1$ and the classical inequality for $|f(z)|$ ($f\in {\mathcal B}$), we have
$$A(z)\leq\frac{r+|a_0|}{1+r|a_0|}+\frac{(1-|a_0|^2)r}{1-r}=1-\frac{(1-|a_0|) A_1(|a_0|,r)}{(1+|a_0|r)(1-r)},
$$
where
$$A_1(|a_0|,r)=(1-|a_0|-|a_0|^2)r^2-(3+|a_0|)r+1.
$$
To prove  the first inequality in Theorem \ref{thm2}, it suffices to prove that $A_1(|a_0|,r)\geq0$  for all $|a_0|\in[0,1)$ and $r\leq r_{a_0}$.
If $|a_0|=\frac{\sqrt{5}-1}{2}$, then we have $1-|a_0|-|a_0|^2=0$ and thus,  $A_1(|a_0|,r)\geq0$ is equivalent to
$$r\leq \frac{2}{3+|a_0|}=\frac{2}{5+\sqrt{5}}=\frac{2}{3+|a_0|+\sqrt{5}(1+|a_0|)}.
$$
For  $|a_0|\in [0,1)\backslash \{\frac{\sqrt{5}-1}{2}\}$, we have $1-|a_0|-|a_0|^2\neq0$ and  thus, we may write 
$$A_1(|a_0|,r)=(1-|a_0|-|a_0|^2)\left(r-\frac{2}{3+|a_0|+\sqrt{5}(1+|a_0|)}\right)\left(r-\frac{2}{3+|a_0|-\sqrt{5}(1+|a_0|)}\right).
$$
Then the desired conclusion can be obtained by a simple analysis on two cases $0\leq |a_0|<\frac{\sqrt{5}-1}{2}$ and $\frac{\sqrt{5}-1}{2}<|a_0|<1$.

To prove the sharpness, we let $a\in[0,1)$ and consider the function
$$ \varphi_a(z)=\frac{a+z}{1+az}=a+(1-a^2)\sum_{n=1}^{\infty}(-a)^{n-1}z^n,\quad z\in\mathbb{D}. 
$$
For this function, we find that
\begin{align*}
 A(r) =\, & \frac{r+a}{1+ra}+\frac{(1-a^2)r}{1-ar}+ \left(\frac{1}{1+a}+\frac{r}{1-r}\right)\frac{(1-a^2)^2r^2}{1-a^2r^2}\\
=\, &1-\frac{1-a}{(1+ar)(1-r)}A_1(a,r),
\end{align*}
where
$$A_1(a,r)=(1-a-a^2)r^2-(3+a)r+1.
$$
The above inequality is bigger than $1$ if and only if $A_1(a,r)<0$.
By the similar analysis in the previous proof, we observe that $A_1(a,r)<0$ if and only if $r>r_a=\frac{2}{3+a+\sqrt{5}(1+a)}$,
which implies the sharpness of the constant $r_{a_0}$ in the first part.

Next we prove the second part. Again, by Lemma \ref{lem4} and the classical inequality for $|f(z)|$, we have
$$ B(z) \leq\left(\frac{r+|a_0|}{1+r|a_0|}\right)^2+\frac{(1-|a_0|^2)r}{1-r}=1-\frac{ (1-|a_0|^2)A_2(|a_0|,r)}{(1+|a_0|r)^2(1-r)},
$$
where
$$A_2(|a_0|,r)=(1-|a_0|^2)r^3-(1+2|a_0|)r^2-2r+1.
$$
To prove  the second inequality in Theorem \ref{thm2}, it suffices to show that $A_2(|a_0|,r)\geq0$ for all $|a_0|\in[0,1)$, only in the case when $0\leq r\leq  r_{a_0}'$.
Elementary calculations show that
$$A_2\left(|a_0|,\frac{1}{3}\right)=\frac{1}{27}(1-|a_0|)(7+|a_0|)>0,
$$
$$A_2\left(|a_0|,\frac{1}{2+|a_0|}\right)=-\frac{(1-|a_0|)(1+|a_0|)^2}{(2+|a_0|)^3}<0
$$
and
$$\frac{\partial A_2}{\partial r}=-3|a_0|^2r^2-4|a_0|r-(1-r)(1+3r)-1<0.
$$
Therefore, the desired conclusion follows easily. The sharpness of the constant $r_{a_0}'$ can be established as in the previous case
and thus, we omit the details. The proof of the theorem is complete.
\epf

\section{Refined Bohr inequalities for symmetric functions}  \label{sec4}

\begin{theorem}  \label{thm3}
Suppose that $f\in {\mathcal B}$ and $f(z)=\sum_{n=0}^{\infty} a_{pn+m} z^{pn+m}$, where $p\in\mathbb{N}$ and $0\leq m\leq p$.
Then the following sharp inequality holds:
$$ C(r):=\sum_{n=0}^{\infty}|a_{pn+m}|r^{pn+m}+ \left(\frac{1}{1+|a_m|}
+\frac{r^p}{1-r^p}\right)\sum_{n=1}^{\infty} |a_{pn+m}|^2 r^{2pn+m}\leq1~\mbox{ for $r\leq r_{p,m,a_m}$},
$$
where $r_{p,m,a_m}$ is the unique positive root of  the equation
$$(1-|a_m|-|a_m|^2)r^{p+m}+r^p+|a_m|r^m-1=0.
$$
Further, we have $r_{p,m,a_m}\geq\sqrt[p]{1/(2+|a_m|)}$.
\end{theorem}


\bpf
Simple observation shows that the function $f$ can be represented as
$$f(z)=z^{m} g\left(z^{p}\right),
$$
where $g\in {\mathcal B}$ and $g(z)=\sum_{n=0}^{\infty} b_{n} z^{n}$ with $b_{n}=a_{pn+m} .$
It follows from Lemma \ref{lem4} with $N=1$ that
\begin{align*}
C(r)= \,& r^m\left[|b_0|+\sum_{n=1}^{\infty}|b_{n}|\left(r^{p}\right)^n+ \left(\frac{1}{1+|b_0|}
+\frac{r^p}{1-r^p}\right)\sum_{n=1}^{\infty} |b_{n}|^2 \left(r^{p}\right)^n\right]\\
\leq \, &r^m\left[|b_0|+(1-|b_0|^2)\frac{r^p}{1-r^p}\right]=r^m\left[|a_m|+(1-|a_m|^2)\frac{r^p}{1-r^p}\right]=:A_3(r)\\
\leq \, &A_3\left(r_{p,m,a_m}\right)=1 \quad \text{for}\quad r\leq r_{p,m,a_m}.
\end{align*}
Since $A_3(r)$ is monotonically increasing in $[0,1)$ and
$$A_3(\sqrt[p]{1/(2+|a_m|)})=(2+|a_m|)^{-m/p}\leq1,
$$
we find that $r_{p,m,a_m}\geq\sqrt[p]{1/(2+|a_m|)}$.

For the sharpness, we consider
$$f(z)=z^m\left (\frac{a-z^p}{1-az^p}\right ) =az^m-(1-a^2)\sum_{n=1}^{\infty}a^{n-1}z^{pn+m},\quad a\in[0,1).
$$
 For this function, we have
\begin{align*}
C(r)=\, & r^m\left[a+(1-a^2)\sum_{n=1}^{\infty}a^{n-1}r^{pn}+ \left(\frac{1}{1+a}
+\frac{r^p}{1-r^p}\right)(1-a^2)^2\sum_{n=1}^{\infty} a^{2(n-1)}r^{2pn}\right]\\
=\, &r^m\left[a+(1-a^2)\frac{r^p}{1-r^p}\right],
\end{align*}
which is bigger than $1$  if and only if $r>r_{p,m,a}$.
This completes the sharpness.
\epf

  If we set $m=0$ in Theorem \ref{thm3} and use its proof, then we have following results which improve \cite[Corollary~1]{KP2018-1}
and generalizes  Theorem \Ref{thmB}.

\begin{cor} \label{cor1}
If $f\in {\mathcal B}$ and $f(z)=\sum_{n=0}^{\infty} a_{pn} z^{pn}$ for some  $p\in\mathbb{N}$, then the following sharp inequalities hold:
\begin{enumerate}
	\item[{\rm (a)}] $\ds \sum_{n=0}^{\infty}|a_{pn}|r^{pn}+ \left(\frac{1}{1+|a_0|}
+\frac{r^p}{1-r^p}\right)\sum_{n=1}^{\infty} |a_{pn}|^2 r^{2pn}\leq1$ \quad \text{for}~~$\ds r\leq\frac{1}{\sqrt[p]{2+|a_0|}}.$
	\item[{\rm (b)}] $\ds |a_0|^2+\sum_{n=1}^{\infty}|a_{pn}|r^{pn}+ \left(\frac{1}{1+|a_0|}
+\frac{r^p}{1-r^p}\right)\sum_{n=1}^{\infty} |a_{pn}|^2 r^{2pn}\leq1$ \quad \text{for}~~$\ds r\leq\frac{1}{\sqrt[p]{2}}.$
\end{enumerate}
\end{cor}

%

\begin{cor} \label{cor2}
Suppose that $f\in {\mathcal B}$ and  $f(z)=\sum_{n=1}^{\infty} a_{pn} z^{pn}$ for some  $p\in\mathbb{N}$. Then the following sharp inequalities hold.
\begin{enumerate}
\item[{\rm (A)}] $\ds \sum_{n=1}^{\infty}|a_{pn}|r^{pn}+ \left(\frac{1}{1+|a_p|}
+\frac{r^p}{1-r^p}\right)\sum_{n=2}^{\infty} |a_{pn}|^2 r^{p(2n-1)}\leq1 \quad \text{for}~~r\leq\sqrt[p]{\frac{2}{\alpha(|a_p|)}},
$
where
$$\alpha(|a_p|)=1+|a_p|+\sqrt{(1-|a_p|)(5+3|a_p|)}.
$$
In particular, this inequality holds for $r\leq\sqrt[p]{3/5}$.
\item[{\rm (B)}] $\ds \sum_{n=1}^{\infty}|a_{pn}|r^{pn}+ \left(\frac{1}{1+|a_p|}
+\frac{r^p}{1-r^p}\right)\sum_{n=1}^{\infty} |a_{pn}|^2 r^{p(2n-1)}\leq1 \quad \text{for}~~r\leq\sqrt[p]{\frac{5-\sqrt{17}}{2}}.
$
\end{enumerate}
\end{cor}
\bpf The inequality (A) follows if we let $m=p$ in Theorem \ref{thm3}. Simple analysis shows that
$$\inf_{|a_p|\in[0,1)}\frac{2}{1+|a_p|+\sqrt{(1-|a_p|)(5+3|a_p|)}}=\frac{3}{5}, ~\mbox{ i.e. }~
 \sqrt[p]{\frac{2}{\alpha(|a_p|)}}\geq \sqrt[p]{\frac{3}{5}}.
$$

Next we prove (B). It follows from the proof of Theorem \ref{thm3} for the case $m=p$ that
\begin{align*}
&\sum_{n=1}^{\infty}|a_{pn}|r^{pn}+ \left(\frac{1}{1+|a_p|}
+\frac{r^p}{1-r^p}\right)\sum_{n=1}^{\infty} |a_{pn}|^2 r^{p(2n-1)}\\
=\,& \sum_{n=1}^{\infty}|a_{pn}|r^{pn}+ \left(\frac{1}{1+|a_p|}
+\frac{r^p}{1-r^p}\right)\sum_{n=2}^{\infty} |a_{pn}|^2 r^{p(2n-1)}+\left(\frac{1}{1+|a_p|}
+\frac{r^p}{1-r^p}\right)|a_p|^2r^{p}\\
\leq\, &r^p\left(|a_p|+(1-|a_p|^2)\frac{r^p}{1-r^p}+\frac{|a_p|^2}{1+|a_p|}
+\frac{|a_p|^2r^p}{1-r^p}\right)\\
=\,&|a_p|r^p+\frac{r^{2p}}{1-r^p}+\frac{|a_p|^2r^p}{1+|a_p|}.
\end{align*}
The remainder of the proof is similar to that of the proof of { \rm(\cite[Theorem~2 (b)]{PVW})} and it follows
just by replacing $|a_1|$ and $r$, respectively, by $|a_p|$ and $r^p$, and the extremal function turned out to be $z^p\left (\frac{a-z^p}{1-az^p}\right )$.
So we omit the details.
\epf

Corollary \ref{cor2}(A) not only generalizes  \cite[Theorem~2 (a)]{PVW}, but also improves upon it.
Furthermore, Corollary \ref{cor2}(B) generalizes  \cite[Theorem~2 (b)]{PVW}.

\section{Refined Bohr type inequalities}  \label{sec5}

The following result is a further refinement corresponding to Theorem \Ref{thmE}.

\begin{theorem} \label{thm4}
Suppose that $f\in {\mathcal B}$, $f(z)=\sum_{n=0}^{\infty} a_{n} z^{n}$ and $S_{r}$ denotes the Riemann surface of the
function $f^{-1}$ defined on the image of the subdisk $|z|<r$ under the mapping $f$.
Then
$$ D(r):=\sum_{n=0}^{\infty}\left|a_{n}\right| r^{n}+\left(\frac{1}{1+|a_0|}
+\frac{r}{1-r}\right)\sum_{n=1}^{\infty} |a_n|^2 r^{2n}+\frac{8}{9}\left(\frac{S_{r}}{\pi}\right) \leq 1 \text { for } r \leq \frac{1}{3}
$$
and the numbers $1/3$ and $8/9$ cannot be improved. Moreover,
$$ E(r):=|a_0|^2+\sum_{n=1}^{\infty}\left|a_{n}\right| r^{n}+\left(\frac{1}{1+|a_0|}
+\frac{r}{1-r}\right)\sum_{n=1}^{\infty} |a_n|^2 r^{2n}+\frac{9}{8}\left(\frac{S_{r}}{\pi}\right) \leq 1 \text { for } r \leq \frac{1}{3-a}
$$
and the constant $9/8$ cannot be improved.
\end{theorem}
\bpf
It follows from Lemma \ref{lem2} that
$$\frac{S_r}{\pi}=\frac{1}{\pi}\int\int_{|z|<r}|f'(z)|^2dxdy=\sum_{n=1}n|a_n|^2r^{2n}\leq\frac{(1-|a_0|^2)^2r^2}{(1-|a_0|^2r^2)^2}\quad \text{for}~~0\leq r\leq
\frac1{\sqrt{2}}.
$$
We now consider the first part. For the case $r\leq1/3$, by Lemma \ref{lem4}, we have
$$
D(r) \leq |a_0|+\frac{(1-|a_0|^2)r}{1-r}+\frac{8(1-|a_0|^2)^2r^2}{9(1-|a_0|^2r^2)^2}=: A_4(r)
$$
so that, because $A_4(r)$ is increasing,
$$D(r) \leq  A_4\left(\frac{1}{3}\right)=1-\frac{(1-|a_0|)^3}{2(9-|a_0|^2)^2}(5+|a_0|)(13+4|a_0|-|a_0|^2)\leq1.
$$
To prove that the constant $8/9$ is sharp, we consider the function
$$f(z)=\phi_a(z)=\frac{a-z}{1-az}.
$$
For this function, straightforward calculations show that
\begin{align*}
&\sum_{n=0}^{\infty}\left|a_{n}\right| r^{n}+\left(\frac{1}{1+|a_0|}
+\frac{r}{1-r}\right)\sum_{n=1}^{\infty} |a_n|^2 r^{2n}+\lambda\left(\frac{S_{r}}{\pi}\right)\\
=\, & a+\frac{(1-a^2)r}{1-ar}+ \frac{1+ar}{(1+a)(1-r)}\frac{(1-a^2)^2r^2}{1-a^2r^2}+\lambda\frac{(1-a^2)^2r^2}{(1-a^2r^2)^2}.
\end{align*}
For $r=1/3$, the right hand side of the last expression becomes
$$1+\frac{(1-a)^2}{2(9-a^2)^2}\{8(9\lambda-8)-8(9\lambda+4)(1-a)+6(3\lambda+2)(1-a)^2+4(1-a)^3-(1-a)^4\},
$$
which is easily seen to be bigger than $1$ in the case $\lambda>8/9$ and $a\rightarrow1$.

Again, if $r\leq1/(3-|a_0|)$, then as in the previous case it follows from Lemma \ref{lem4} that
\begin{align*}
E(r)~\leq& ~|a_0|^2+\frac{(1-|a_0|^2)r}{1-r}+\frac{9(1-|a_0|^2)^2r^2}{8(1-|a_0|^2r^2)^2}=: A_5(r)\\
~\leq& ~ A_5\left(\frac{1}{3-|a_0|}\right)=1-\frac{(1-|a_0|)^3(1+|a_0|)}{8(2-|a_0|)(3-2|a_0|)^2}[(54-39|a_0|)+|a_0|^2(6-|a_0|)]\\
~\leq& ~1.
\end{align*}
Finally, to prove that the constant $9/8$ is sharp in the second part of the statement of the theorem, we consider the function $f(z)=\phi_a(z)$
which is given above. Again, for this function, straightforward calculations show that
\begin{align*}
&|a_0|^2+\sum_{n=1}^{\infty}\left|a_{n}\right| r^{n}+\left(\frac{1}{1+|a_0|}
+\frac{r}{1-r}\right)\sum_{n=1}^{\infty} |a_n|^2 r^{2n}+\lambda\left(\frac{S_{r}}{\pi}\right)\\
=\, & a^2+\frac{(1-a^2)r}{1-ar}+ \frac{1+ar}{(1+a)(1-r)}\frac{(1-a^2)^2r^2}{1-a^2r^2}+\lambda\frac{(1-a^2)^2r^2}{(1-a^2r^2)^2}.
\end{align*}
For $r=1/(3-a)$, the last expression reduces to
$$1+\frac{(1-a)^2(1+a)}{9(3-2a)^2(2-a)} A_6(a,\lambda),
$$
which is again seen to be greater than $1$ when $\lambda>9/8$ and $a\rightarrow1$, where
$$A_6(a,\lambda)=(8\lambda-9)+12(\lambda-3)(1-a)+2(\lambda-18)(1-a)^2-3\lambda(1-a)^3-\lambda(1-a)^4.
$$
This completes the proof.
\epf

\begin{rems}
From the proof of Theorem \ref{thm4}, it is easy to see that the number $1/(3-a)$ in Theorem \ref{thm4} is not sharp for all $a\in[0,1)$.
However, it is sharp in the sense of approximation based on the following observation. On one hand, the number
$1/(3-a)$ approaches $1/2$ as $a\rightarrow1$.
On the other hand, suppose that $f\in {\mathcal B}$ and $f(z)=\sum_{n=0}^{\infty} a_{n} z^{n}$, we cannot find any $c>0$ such that
$$|a_0|^2+\sum_{n=1}^{\infty}\left|a_{n}\right| r^{n}+\left(\frac{1}{1+|a_0|}
+\frac{r}{1-r}\right)\sum_{n=1}^{\infty} |a_n|^2 r^{2n}+c\left(\frac{S_{r}}{\pi}\right) \leq 1 \text { for } r \leq \frac{1}{2};
$$
for instance, consider $f(z)=z$ and $r=1/2$.
\end{rems}

 Next, we consider another refined Bohr type inequality  with the square of the modulus of $f(z)-a_0(f)$, instead of area term.

\begin{theorem} \label{thm5}
 Suppose that $f\in {\mathcal B}$ and $f(z)=\sum_{n=0}^{\infty} a_{n} z^{n}$.  Then
$$ F(z):=\sum_{n=0}^{\infty}\left|a_{n}\right| r^{n}+\left(\frac{1}{1+|a_0|}
+\frac{r}{1-r}\right)\sum_{n=1}^{\infty} |a_n|^2 r^{2n}+|f(z)-a_0|^2\leq 1 \text { for } |z|=r \leq \frac{1}{3}
$$
if and only if $0\leq|a_0|\leq4\sqrt{2}-5\approx 0.656854$.
\end{theorem}
\bpf
Assume that $0\leq|a_0|\leq4\sqrt{2}-5$. Then $|a_0|^2+10|a_0|-7\leq0$.
If $4\sqrt{2}-5\geq|a_0|\geq r$ and $r\leq1/3$, then it follows from Lemma \ref{lem1} and Lemma \ref{lem4} with $N=1$ that
\begin{align*}
F(z)~\leq& ~|a_0|+\frac{(1-|a_0|^2)r}{1-r}+\left(\frac{(1-|a_0|^2)r}{1-|a_0|r}\right)^2=: A_7(r)\\
~\leq&~A_7\left(\frac{1}{3}\right)=1+\frac{(1-|a_0|)^2(|a_0|^2+10|a_0|-7)}{2(3-|a_0|)^2}\leq1.
\end{align*}
Again, if $0\leq|a_0|<r\leq1/3$, then it follows from Lemma \ref{lem1} and Lemma \ref{lem4} with $N=1$ that
\begin{align*}
F(z)~\leq& ~|a_0|+\frac{(1-|a_0|^2)r}{1-r}+\left(\frac{\sqrt{1-|a_0|^2}r}{\sqrt{1-r^2}}\right)^2=: A_8(r)\\
~\leq&~A_8\left(\frac{1}{3}\right)=|a_0|+\frac{1}{2}(1-|a_0|^2)+\frac{\sqrt{1-|a_0|^2}}{8}\\
~\leq&~\frac{1}{3}+\frac{1}{2}+\frac{1}{8}<1.
\end{align*}

To complete the proof, we choose $f(z)=\phi_a(z)=(a-z)/(1-az)$ and for this function, we have that
\begin{align*}
F(r)~=&~a+\frac{(1-a^2)r}{1-r}+\left(\frac{(1-a^2)r}{1-ar}\right)^2
\end{align*}
which for $r=1/3$ reduces to
$$1+\frac{(1-a)^2(a^2+10a-7)}{2(3-a)^2},
$$
and this is larger than $1$ if and only if $4\sqrt{2}-5<a<1$.
\epf

Clearly, Theorem \ref{thm5} improves \cite[Theorem~3]{KP2018-2} partly.
Surprisingly, if the square term $|f(z)-a_0|^2$ in Theorem \ref{thm5} is replaced by $|f(z)-a_0|$, then
we have the following result which works for all $|a_0|\in[0,1)$.

\begin{theorem} \label{thm6}
 Suppose that $f\in {\mathcal B}$ and $f(z)=\sum_{n=0}^{\infty} a_{n} z^{n}$.  Then
$$ G(z):=\sum_{n=0}^{\infty}\left|a_{n}\right| r^{n}+\left(\frac{1}{1+|a_0|}
+\frac{r}{1-r}\right)\sum_{n=1}^{\infty} |a_n|^2 r^{2n}+|f(z)-a_0|\leq 1 \text { for } |z|=r \leq \frac{1}{5}
$$
and the number $1/5$  cannot be improved. Moreover,
$$ H(z):=|a_0|^2+\sum_{n=1}^{\infty}\left|a_{n}\right| r^{n}+\left(\frac{1}{1+|a_0|}
+\frac{r}{1-r}\right)\sum_{n=1}^{\infty} |a_n|^2 r^{2n}+|f(z)-a_0|\leq 1
$$
$\text { for } |z|=r \leq \frac{1}{3}$,
and the constant $1/3$ cannot be improved.
\end{theorem}
\bpf
For the first part, we first consider $|a_0|\geq r$ and $r\leq1/5$. Then it follows from Lemma \ref{lem1} and  Lemma \ref{lem4} with $N=1$ that
\begin{align*}
G(z)~\leq&~|a_0|+\frac{(1-|a_0|^2)r}{1-r}+\frac{(1-|a_0|^2)r}{1-|a_0|r}=: A_9(r)\\
~\leq&~A_9\left(\frac{1}{5}\right)=1-\frac{(1-|a_0|)^2(11-|a_0|)}{4(5-|a_0|)}\leq1.
\end{align*}
If $|a_0|<r\leq1/5$, then it follows from Lemmas \ref{lem1} and   \ref{lem4} with $N=1$ that
\begin{align*}
G(z)~\leq&~|a_0|+\frac{(1-|a_0|^2)r}{1-r}+\frac{\sqrt{1-|a_0|^2}r}{\sqrt{1-r^2}}=: A_{10}(r)\\
~\leq&~A_{10}\left(\frac{1}{5}\right)=|a_0|+\frac{1}{4}(1-|a_0|^2)+\frac{\sqrt{1-|a_0|^2}}{2\sqrt{6}}\\
~\leq&~\frac{1}{5}+\frac{1}{4}+\frac{1}{4}<1.
\end{align*}
To check its sharpness, we choose $f(z)=\phi_a(z)=(a-z)/(1-az)$ and as before for this function, we find that
\begin{align*}
G(r)~=&~a+\frac{(1-a^2)r}{1-r}+\frac{(1-a^2)r}{1-ar}\\
~=&~1-\frac{(1-a)}{(1-r) (1-ar)}[r^2a^2+(3 r^2-3 r)a+(r^2-3 r+1)]
\end{align*}
which is larger than $1$ if and only if
$$ A_{11}(a,r):=r^2a^2+(3 r^2-3 r)a+(r^2-3 r+1)<0.
$$
If $1/5< r <(3-\sqrt{5})/2$, we let $a_r=\frac{3(1-r)-\sqrt{5r^2-6 r+5}}{2r}$.
Elementary calculations show that  $a_r\in(0,1)$ and $\frac{3r-3r^2}{2r^2}>1$, and thus
$A_{11}(a,r)<A_{11}(a_r,r)=0$ when $a_r<a<1$. This proves the sharpness.

For the proof of the second part of the theorem,  as in the previous case, we have
\begin{align*}
H(z)~\leq&~|a_0|^2+\frac{(1-|a_0|^2)r}{1-r}+\frac{(1-|a_0|^2)r}{1-|a_0|r}=: A_{12}(r)\\
~\leq&~A_{12}\left(\frac{1}{3}\right)=1-\frac{(1-|a_0|)^2(1+|a_0|)}{2(3-|a_0|)}\leq1
\end{align*}
when $|a_0|\geq r$ and $r\leq1/3$; and
\begin{align*}
H(z)~\leq&~|a_0|^2+\frac{(1-|a_0|^2)r}{1-r}+\frac{\sqrt{1-|a_0|^2}r}{\sqrt{1-r^2}}=: A_{13}(r)\\
\leq&A_{13}\left(\frac{1}{3}\right)=|a_0|^2+\frac{1}{2}(1-|a_0|^2)+\frac{\sqrt{2}}{4}\sqrt{1-|a_0|^2}\\
\leq&\frac{1}{9}+\frac{1}{2}+\frac{\sqrt{2}}{4}<1
\end{align*}
when $|a_0|<r\leq1/3$. To check its sharpness, we choose again $f(z)=\phi_a(z)=(a-z)/(1-az)$ and for this function we find that
\begin{align*}
H(r)~=&~a^2+\frac{(1-a^2)r}{1-r}+\frac{(1-a^2)r}{1-ar}\\
~=&~1-\frac{(1-a^2)}{(1-r) (1-ar)}[(1-3r+r^2)+(-r+2r^2)a]
\end{align*}
which is larger than $1$ if and only if
$ A_{14}(a,r):=(1-3r+r^2)+(-r+2r^2)a<0.
$
If $1/3< r <(3-\sqrt{5})/2$, then we let $a_r=\frac{r^2-3 r+1}{r-2 r^2}$ and we see that
$a_r\in(0,1)$ and $A_{14}(a,r)<A_{14}(a_r,r)=0$ when $a_r<a<1$.
This completes the sharpness.
\epf

Finally, we present an improved version of Theorem \Ref{thmF}.

\begin{theorem} \label{thm7}
Suppose that $f\in {\mathcal B}$ and $f(z)=\sum_{n=0}^{\infty} a_{n} z^{n}$.
Then
$$ I(z):=|f(z)|+\left|f'(z)\right|r+\sum_{n=2}^{\infty}|a_n|r^n+ \left(\frac{1}{1+|a_0|}
+\frac{r}{1-r}\right)\sum_{n=1}^{\infty} |a_n|^2 r^{2n} \leq 1
$$
for $|z|=r\leq\frac{\sqrt{17}-3}{4}$ and the constant $\frac{\sqrt{17}-3}{4}$ is best possible.
Moreover,
$$ J(z):=|f(z)|^2+\left|f'(z)\right|r+\sum_{n=2}^{\infty}|a_n|r^n+ \left(\frac{1}{1+|a_0|}
+\frac{r}{1-r}\right)\sum_{n=1}^{\infty} |a_n|^2 r^{2n} \leq 1
$$
for $|z|=r\leq r_0,$
where $r_0\approx0.385795$ is the unique positive root of the equation
$$1-2r-r^2-r^3-r^4=0
$$ and $r_0 $ is best possible.
\end{theorem}
\bpf
By simple calculation we know that $\frac{2r}{1-r^2}\leq1$ if $0\leq r\leq\sqrt{2}-1$.
Combining Lemma \ref{lem4}, the classical inequality for $|f(z)|$ and the Schwarz-Pick lemma, we have
\begin{align*}
I(z)~\leq&~|f(z)|+\frac{|z|}{1-|z|^2}(1-|f(z)|^2)+\frac{(1-|a_0|^2)r^2}{1-r}\\
~\leq&~ \frac{r+|a_0|}{1+r|a_0|}+\frac{r}{1-r^2}\left[1-\left(\frac{r+|a_0|}{1+r|a_0|}\right)^2\right]+\frac{(1-|a_0|^2)r^2}{1-r}\\
~=&~\frac{r+|a_0|}{1+r|a_0|}+\frac{r(1-|a_0|^2)}{(1+r|a_0|)^2} +\frac{(1-|a_0|^2)r^2}{1-r}\\
~=&~1+\frac{1-|a_0|}{(1+|a_0|r)^2(1-r)}[-1+3r-r^2+(2r^2+r^3)|a_0|+(2r^3+r^4)|a_0|^2+r^4|a_0|^3]\\
~\leq&~1+\frac{1-|a_0|}{(1+|a_0|r)^2(1-r)}[-1+3r-r^2+(2r^2+r^3)+(2r^3+r^4)+r^4]\\
~=&~1+\frac{2(1-|a_0|)(1+r^2)}{(1+|a_0|r)^2(1-r)}\left(r-\frac{\sqrt{17}-3}{4}\right)\left(r+\frac{\sqrt{17}+3}{4}\right)
\end{align*}
which is no more than $1$ for $0\leq r\leq \frac{\sqrt{17}-3}{4}<\sqrt{2}-1$. It is worth pointing out in the second inequality above, we have used the fact
that
$$\Phi(X)=X+\lambda (1-X^2)\leq \Phi(X_0) ~\mbox{ when $\ds X=|f(z)|\leq X_0=\frac{r+|a_0|}{1+r|a_0|}$ and $\ds \lambda =\frac{|z|}{1-|z|^2}$.}
$$
Clearly, $1-\frac{r}{1-r^2}\geq0$ if $0\leq r\leq\frac{1}{2}$.
Again, as in the previous case, we have
\begin{align*}
J(z)
~\leq&~ \left(\frac{r+|a_0|}{1+r|a_0|}\right)^2+\frac{r(1-|a_0|^2)}{(1+r|a_0|)^2}+\frac{(1-|a_0|^2)r^2}{1-r}\\
~=&~ 1-\frac{(1-r^2)(1-|a_0|^2)}{(1+r|a_0|)^2}+\frac{r(1-|a_0|^2)}{(1+r|a_0|)^2}+\frac{(1-|a_0|^2)r^2}{1-r}\\
~=&~1+\frac{1-|a_0|^2}{(1+|a_0|r)^2(1-r)}[-(1-r)(1-r-r^2)+r^2(1+r|a_0|)^2]\\
~\leq&~1+\frac{1-|a_0|^2}{(1+|a_0|r)^2(1-r)}[-1+2r+r^2+r^3+r^4]\leq1
\end{align*}
for $0\leq r\leq r_0<1/2$.

To check its sharpness, by choosing $f(z)=\varphi_a(z)=(a+z)/(1+az)$, we get
\begin{align*}
J(r)=\, & \left(\frac{r+a}{1+ra}\right)^2+\frac{(1-a^2)r}{(1+ar)^2}+\frac{(1-a^2)ar^2}{1-ar}+ \frac{1+ar}{(1+a)(1-r)}\frac{(1-a^2)^2r^2}{1-a^2r^2}\\
=\, &1+\frac{1-a^2}{(1+ar)^2(1-r)}[-1+2r+r^2-r^3+2r^3a+r^4a^2].
\end{align*}
The last expression is larger than 1 if and only if
$$A_{15}(a,r):=-1+2r+r^2-r^3+2r^3a+r^4a^2>0
$$
for all $a\in[0,1)$ and all ~$r$~¡¡in some subset of $[0,1)$.
The equivalent condition implies $A_{15}(a,r)\geq0$ by $a\rightarrow1$.
It is easy to see that $A_{15}(1,r)$ is  greater than of equal to zero if and only if $r\geq r_0$.
This shows the sharpness.
\epf

\subsection*{Conflict of Interests}
The authors declare that there is no conflict of interests regarding the publication of this paper.

\subsection*{Acknowledgments}
The research of the first author was supported by Application-Oriented Characterized Disciplines,
Double First-Class University Project of Hunan Province (Xiangjiaotong [2018]469) and
the Science and Technology Plan Project of Hunan Province (No. 2016TP1020).
The second author was supported by the NSFs of China~(No. 11961013) 
and the Natural Science Foundation of Guangxi (No. 2018GXNS-FAA050005).
The  work of the third author is supported by Mathematical Research Impact Centric Support (MATRICS) of
the Department of Science and Technology (DST), India  (MTR/2017/000367).

\end{document}